\def\draft{n}
\documentclass{amsart}
\usepackage{fullpage,amssymb,epic,eepic,epsfig,amscd,pb-diagram,lamsarrow,pb-lams}

\theoremstyle{plain}

\newtheorem{theorem}{Theorem}
\newtheorem{proposition}{Proposition}[section]
\newtheorem{lemma}[proposition]{Lemma}

\theoremstyle{definition}
\newtheorem{definition}[proposition]{Definition}

\newtheorem{question}{Question}

\theoremstyle{remark}

\newtheorem{remark}[proposition]{Remark}

\def\printname#1{
	\if\draft y
		\smash{\makebox[0pt]{\hspace{-0.5in}
			\raisebox{8pt}{\tt\tiny #1}}}
	\fi
}

\newcommand{\psdraw}[2]
         {\begin{array}{c} \hspace{-1.3mm}
	\raisebox{-4pt}{\epsfig{figure=draws/#1.eps,width=#2}}
	\hspace{-1.9mm}\end{array}}

\newlength{\standardunitlength}
\catcode`\@=11
\long\def\@makecaption#1#2{%
     \vskip 10pt

\setbox\@tempboxa\hbox{%\ifvoid\tinybox\else\box\tinybox\fi
       \small\sf{\bfcaptionfont #1. }\ignorespaces #2}%
     \ifdim \wd\@tempboxa >\captionwidth {%
         \rightskip=\@captionmargin\leftskip=\@captionmargin
         \unhbox\@tempboxa\par}%
       \else
         \hbox to\hsize{\hfil\box\@tempboxa\hfil}%
     \fi}
\font\bfcaptionfont=cmssbx10 scaled \magstephalf
\newdimen\@captionmargin\@captionmargin=2\parindent
\newdimen\captionwidth\captionwidth=\hsize
\catcode`\@=12

\def\lbl#1{\label{#1}\printname{#1}}

\def\BZ{\mathbb Z}
\def\BQ{\mathbb Q}

\def\A{\mathcal A}

\def\O{\mathcal O}

\def\l{\lambda}
\def\Ga{\Gamma}

\def\fti{finite type invariant}

\def\Ga{\Gamma}

\def\b{\beta}

\def\sub{\subset}

\def\lk{{\text{lk}}}

\def\ti{\widetilde}

\def\st#1#2#3{\overset{#1}{\underset{#2}{
\begin{array}{c} \hspace{-1.3mm}
	\raisebox{-4pt}{\psfig{figure=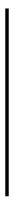,height=0.2in} }
	\hspace{-1.9mm}\end{array} }}#3}

\def\AS{\mathrm{AS}}
\def\IHX{\mathrm{IHX}}

\def\longto{\longrightarrow}
\def\bmu{\bar{\mu}}
\def\om{\omega}
\def\Ztr{Z^{\mathrm{tr}}}
\def\Zmin{Z^{\mathrm{min}}}
\def\At{\A^{\mathrm{tr}}}

\begin{document}

%%%%%%%%%%%%%%%%%%%%%%{page1}

\title[Links with trivial Alexander module and nontrivial Milnor 
invariants]{Links with trivial Alexander module and nontrivial Milnor 
invariants}

\author{Stavros Garoufalidis}
\address{School of Mathematics \\
          Georgia Institute of Technology \\
          Atlanta, GA 30332-0160, USA. }
\email{stavros@math.gatech.edu}

\thanks{Supported by NSF grant DMS-02-03129 and an Israel-US BSF grant. \\
        This and related preprints can also be obtained at
{\tt http://www.math.gatech.edu/$\sim$stavros } 
%and {\tt http: }
\newline
1991 {\em Mathematics Classification.} Primary 57N10. Secondary 57M25.
}

\date{
This edition: June 7, 2002 \hspace{0.5cm} First edition: June 7, 2002.}

%\dedicatory{Preliminary notes. Please do not distribute under
%any circumstances!}

\begin{abstract}
Cochran constructed many links with Alexander module that of the unlink
and some nonvanishing Milnor invariants, using as input commutators in
a free group and as an invariant the longitudes of the links. We present a 
different and conjecturally complete construction, that uses elementary 
properties of clasper surgery, and a different invariant, the tree-part of 
the LMO invariant. Our method also constructs links with trivial higher
Alexander modules and nontrivial Milnor invariants.
\end{abstract}

\maketitle

%\tableofcontents

%%%%%%%%%%%%%%%%% the text file

\section{Introduction}
\lbl{sec.intro}
\subsection{History of the problem}
\lbl{sub.history}

Two of the best studied topological invariants a link $L$ in $S^3$ are its 
{\em 
Alexander module} $A(L)$ which measures the homology of the universal abelian
cover of $S^3-L$, and its collection of {\em Milnor invariants} 
$\bmu(L)$, which are concordance (and sometimes link homotopy) 
invariants, defined
modulo a recursive indeterminacy. Let us say that $L$ has {\em trivial} 
Alexander module (resp. Milnor invariants) if $A(L)=A(\O)$ 
(resp. $\bmu(L)=\bmu(\O)=0$) for an unlink $\O$. Despite the indeterminacy
of the Milnor invariants, note that the vanishing of all Milnor
invariants is a well-defined statement.

Using the language of {\em longitudes} $\l_i$ of components of $L$,
Milnor showed that a link $L$ has vanishing Milnor invariants iff
$\l_i(L) \subset \pi_{\om}$ for all $i$, where $\pi=\pi_1(S^3-L)$ and 
$\pi_{\om}=\cap_{n=1}^\infty \pi_n$ is the intersection of the {\em lower 
central series} $\pi_n$ of $\pi$, defined by $\pi_1=\pi$ and 
$\pi_{n+1}=[\pi_n,\pi]$, see \cite{Mi}. $L$ has trivial Alexander module iff 
there is a map $\pi \to F/[[F,F],[F,F]]$ which induces an isomorphism
$\pi/[[\pi,\pi],[\pi,\pi]] \cong F/[[F,F],[F,F]]$.
%Further, the longitudes of a link with trivial Alexander module
%lie in $[[\pi,\pi],[\pi,\pi]]$, \cite[p. 75]{Hi}.

It is natural to ask how independent are the conditions of trivial Alexander
module and trivial Milnor invariants. In a sense, this question asks for
a comparison between the lower central series and the commutator series
of a link group.

In one direction, Levine showed that the vanishing of the Milnor invariants 
of a link $L$ implies that a localization $A(L)_S$ of its Alexander module
(although not the Alexander module itself) vanishes, where $S \subset
\BZ[t_1^{\pm 1}, \dots, t_r^{\pm 1}]$ is the multiplicative set of 
polynomials that evaluate
to $\pm 1$ at $t_1=\dots=t_r=1$; see \cite{Le}.
A boundary link has vanishing Milnor invariants, and its Alexander module
splits as a direct sum of a trivial module and a torsion module. It was shown
in \cite{GL2} that all torsion modules with the appropriate symmetry can
be realized.

In the opposite direction, if $L$ has trivial Alexander module, then it
is known that some low order Milnor invariants vanish, \cite{Le,Tr}. For
example, all nonrepeated (link homotopy) invariants with at most $5$ indices
vanish. On the other hand, Cochran  
constructed a class of links with trivial Alexander module
and nontrivial Milnor invariants; such links are not even be concordant
to homology boundary links.

Cochran's construction used iteration, and used as a pattern certain
elements in the lower central series of the free group. There is enough
explicitness and control on the iteration that enabled Cochran to compute
the longitudes directly and verify that these links have vanishing Alexander
modules. Further, a geometric interpretation of Milnor invariants in terms
of cycles on Seifert surfaces allowed Cochran to conclude that the 
constructed links have nontrivial Milnor invariants.

As an elementary application of the calculus of claspers, we will construct
a plethora of links with vanishing Alexander module. For these links, we can
compute the tree part of the LMO invariant (which can be identified with
Milnor invariants, \cite{HM}), using formal Gaussian integration.
As a result, we will construct many (and conjecturally all) links with
trivial Alexander module and nontrivial Milnor invariants.
The next definition explains the patterns that we will use in our
construction.

\begin{definition}
\lbl{def.pattern}
Let $\At(r)$ (or simply, $\At$, in case $r$ is clear) denote the vector space 
over $\BQ$ generated by vertex-oriented
unitrivalent trees, whose univalent vertices are labeled by $r$ colors,
modulo the $\AS$ and $\IHX$ relation. $\At(r)$ is a graded vector space,
where the degree of a graph is half the number of vertices.
We will call a tree of degree $1$ 
(with two univalent vertices and no trivalent ones) a {\em strut}.

A {\em pattern} $\b$ is an element of $\At(r)$ which is represented by a 
tree which has a trivalent vertex $v$ such that $\b-v$ has no strut 
components.
\end{definition}

The next figure gives some examples of nonvanishing patterns:
$$
\psdraw{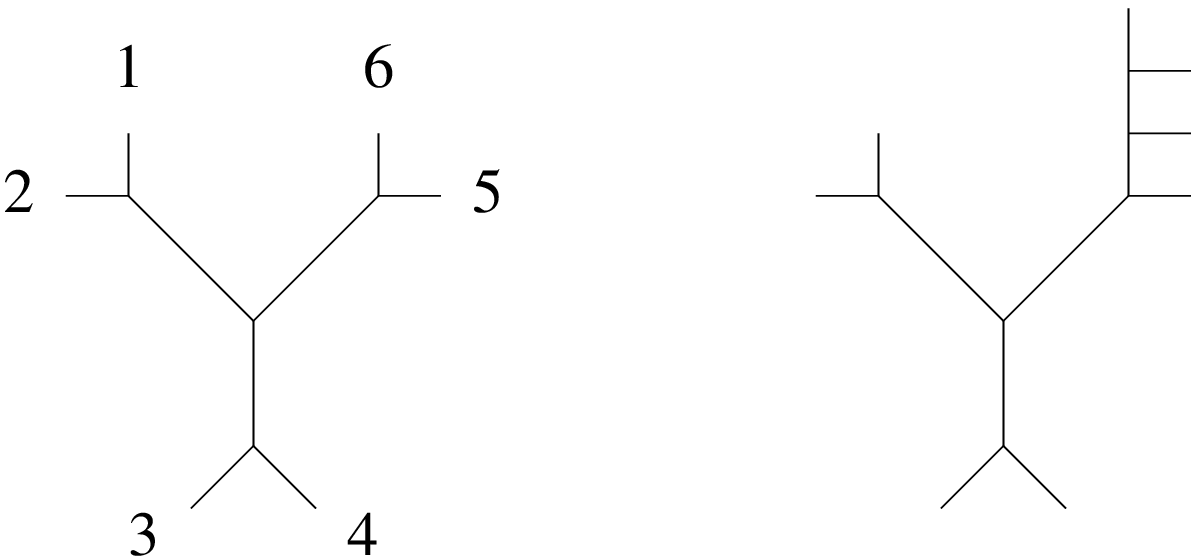}{2in}
$$

\begin{theorem}
\lbl{thm.1}
For every nonvanishing pattern $\b \in \At_m(r)$ there exists a link 
$L(\b)$ with $r$ components such that $A(L(\b))=A(\O)$, all Milnor invariants 
of degree less than $m$ vanish and some Milnor invariant of degree $m$ do not.
\end{theorem}

Our construction adapts without change to the case of links with trivial
{\em higher Alexander modules}. Although classical, these modules appeared
only recently in work of Cochran-Orr-Teichner \cite{COT} and subsequent work
of Cochran, \cite{Co2}. Given a group $\pi$, consider its {\em commutator 
series} defined by $\pi^{(0)}=\pi$ and $\pi^{(n+1)}=[\pi^{(n)},\pi^{(n)}]$.

\begin{definition}
\lbl{def.nalex}
We will say that a link $L$ in a homology sphere $M$ has {\em trivial $n$th
Alexander module} if it has a map $\pi \longto F/F^{(n+1)}$
which induces an isomorphism $\pi/\pi^{(n+1)} \cong F/F^{(n+1)}$, where 
$\pi=\pi_1(M-L)$. 
\end{definition}

The next definition explains the $n$-patterns which we will use.

\begin{definition}
\lbl{def.patternn}
Let $c^{(n)}$ be a unitrivalent tree defined by
$$
\psdraw{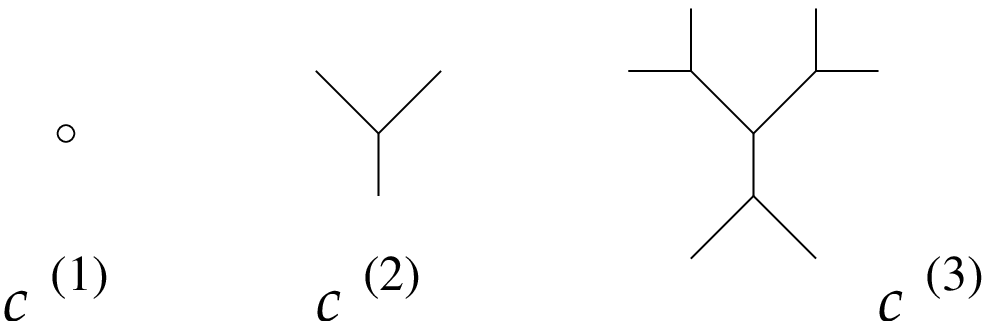}{2in}
$$
In other words, we are adding two univalent vertices in $c^{(n+1)}$ to each 
of the univalent vertices of $c^{(n)}$.
An {\em $n$-pattern} $\b^{(n)}$ is an element of $\At(r)$ which is represented
by a tree $\b^{(n)}$ such that $c^{(n)} \sub \b^{(n)}$ and $\b^{(n)}-c^{(n)}$
has no strut components.
\end{definition}

The proof of Theorem \ref{thm.1} generalizes without change to the following

\begin{theorem}
\lbl{thm.2}
For every nonvanishing $n$-pattern $\b^{(n)} \in \At_m(r)$ there exists a link 
$L(\b^{(n)})$ with $r$ components with trivial $n$th Alexander module, such 
that all Milnor invariants of degree less than $m$ vanish and some Milnor 
invariant of degree $m$ do not.
\end{theorem}

\section{Constucting links by surgery on claspers}
\lbl{sub.construct}

\subsection{What is surgery on a clasper?}
\lbl{sub.surgery}

As we mentioned in the introduction, we will construct links of  
Theorem \ref{thm.2} using {\em surgery on claspers}.
Since claspers play a key role in geometric constructions, as well as
in the theory of \fti s, we include a brief discussion here.
For a reference on claspers and their associated surgery, we refer the reader 
to \cite{Gu2,H} and also to \cite[Section 2]{GGP} (where claspers were called
clovers instead).
It suffices  to say that a clasper is a thickening of a trivalent
graph, and it has a preferred set of loops, called the leaves. The degree of
a clasper is the number of trivalent vertices (excluding those at the 
leaves). With our conventions, the smallest clasper is a
Y-clasper (which has degree one and three leaves), so we explicitly
exclude struts (which would be of degree zero with two leaves).

A clasper $G$ of degree~1 is an embedding $G: N \to M$ of a regular 
neighborhood of the graph $\Ga$ 

$$
\psdraw{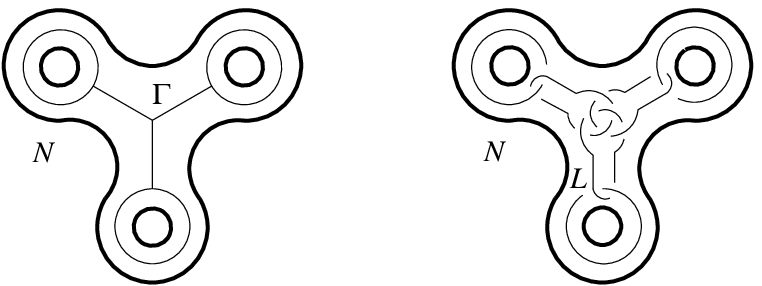}{3in}
$$
in a 3-manifold $M$. Surgery on $G$ can be described by cutting 
$G(N)$ from $M$ (which is a genus~3 handlebody), twisting by a 
fixed diffeomorphism of its boundary (which acts trivially on the homology 
of the boundary) and gluing back. We will denote the result of surgery
by $M_G$.
Alternatively, we can describe surgery on $G$ by surgery on a framed
six component link (the image of $L$) in $M$. The six component link
consists of a $0$-framed Borromean ring and an arbitrarily framed three
component link, the so-called {\em leaves} of $G$.
If one of the leaves bounds a $0$-framed disk disjoint from the rest of $G$,
then surgery on $G$ does not change the ambient 3-manifold $M$, although it
can change an embedded link in $M$. In particular, surgery on a 
clasper of degree~1 is shown as follows:

$$
\psdraw{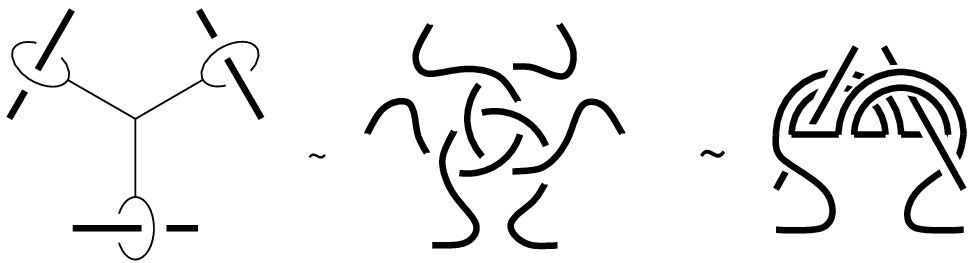}{3.5in}
$$
In general, surgery on a clasper $G$ of degree $n$ can be described in terms
of simultaneous surgery on $n$ claspers $G_1, \dots, G_n$, which are obtained
from $G$ after breaking its edges and inserting Hopf links as follows:

$$
\psdraw{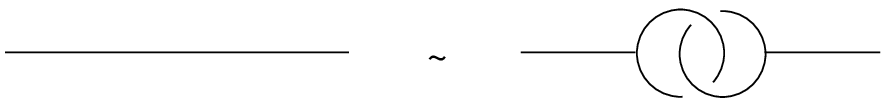}{2in}
$$

\subsection{A basic principle}
\lbl{sub.basic}

Surgery on a clasper is described by twisting by a surface diffeomorphism
that acts trivially on homology, thus we have the basic principle:

\begin{center}
\fbox{Clasper surgery preserves the homology}
\end{center}

Surgery on claspers with leaves of a resticted type has already been studied 
and used successfully in \cite{GR} (where the leaves were assumed null 
homologous in a knot complement), \cite{GL1} (and where the leaves where
null homotopic) and \cite{GK} (where the leaves where in the kernel of a map
to a free group). It is important to study not only 3-manifolds
but rather pairs of 3-manifolds together with a representation of their 
fundamental group into a fixed group. Claspers adapt well to this point
of view, as we explain next.

Consider a pair $(N,\rho)$ of a 3-manifold $N$ (possibly noncompact) and a 
representation $\rho: \pi_1(N) \to \Ga$ for some group $\Ga$. Consider
a clasper $G \subset N$ whose leaves are mapped to $1$ under $\rho$. We will
call such claspers $\rho$-{\em null}, or simply {\em null}, if $\rho$ is
clear. Surgery on $G$ gives rise to a 4-manifold $W$ whose boundary consists
of one copy of $N$ and one copy of $N_G$. We may think that $W$ is obtained 
by attaching
$6n$ 2-handles on $N \times I$, where $n=\text{degree}(G)$. Since the cores 
of these handles lie in the kernel of $\rho$, it follows that $\rho$ extends
over $W$, and in particular restricts to a representation $\rho_G$ on
the end $N_G$ of $W$.

\begin{lemma}
\lbl{lem.rho}
We have $H_\ast(N, \rho) \cong H_\ast(N_G, \rho_G)$.
\end{lemma}

\begin{proof}
Let $\ti N$ (resp. $\ti {N_G}$) denote the cover of $N$ (resp. $N_G$) 
corresponding to $\rho$ (resp. $\rho_G$). Surgery on $G$ is equivalent
to surgery on a collection $\{G_1,\dots,G_k\}$ of degree $1$ claspers,
constructed by inserting Hopf links in the edges of $G$. Each $G_i$
lifts to a collection $\ti G_i$ of claspers in $\ti N$; 
let $\ti G=\ti G_1 \cup \dots \ti G_k$. Then, $\ti {N_G}$ can be identified 
with $(\ti N)_{\ti G}$. Since 
clasper surgery preserves homology, the result follows. 
\end{proof}

We will adapt the above lemma in the following situation.
Suppose that $G$ is a clasper in the complement of an unlink 
$X_0=S^3-\O$ of $r$ components whose leaves are null homologous in $X_0$, and 
let $(M,L)$ denote the result of surgery along $G$ on the pair $(S^3,\O)$. It 
follows that $G$ lifts to a family $\ti G$ of claspers in $\ti X_0$ 
(the universal abelian cover of $X$) and that $\ti X$ is obtained from 
$\ti X_{0}$, by surgery on $\ti G$, where $X=M-L$.  Since 
$A(L)=H_1(\ti X, \ti x)$, and clasper surgery preserves homology, it follows 
that $A(M,L)=A(\O)$.

\begin{remark}
There are two known cases where surgery on a null clasper $G \subset X_0$ gives
rise to a link $(M,L)$ with vanishing Milnor invariants.

\noindent
(a) If the leaves of $G$ are null homotopic in $X_0$, then the 
constructed links would be boundary links, as was observed and used 
in \cite{GK}. Boundary links have vanishing Milnor invariants.

\noindent
(b) If $G$ is a connected clasper with at least one loop, then $(M,L)$
is concordant to $(S^3,\O)$, \cite{GL3} and also \cite{CT}. Concordance
preserves Milnor invariants.
\end{remark}

With a bit more effort, we can arrange that $M=S^3$. For this, it suffices
to assume that each connected component $G_i$ of $G$ has a 0-framed 
leaf $l_i$, such that the union of the leaves $\{l_i\}$ is an unlink
in $S^3$.

To finalize the construction of Theorem \ref{thm.1}, consider a pattern $\b$,
and a vertex $v$ of $\b$ such that $\b-v=T_1 \cup T_2 \cup T_3$
where $T_i$ are rooted trees which are not struts. Each rooted tree $T$
corresponds to an element $\phi(T) \in F$ via a map defined in pictures
by:
$$
\psdraw{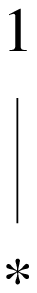}{0.07in} \longto t_1 \in F,
\psdraw{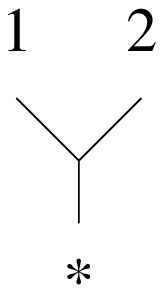}{0.35in} \longto [t_1,t_2] \in F, \hspace{1cm} 
\psdraw{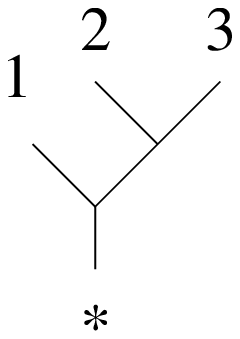}{0.5in}
\longto [t_1, [t_2,t_3]] \in F
$$
If $T$ is not a strut, then $\phi(T) \in [F,F]$. Given $\b$ as above,
we will choose a clasper $G(\b)$ of degree $1$ such that its three leaves
$l_i$ satisfy $l_i=\phi(T_i) \in [F,F]$, for $i=1,2,3$. Then, $L(\b)$
is obtained from the unlink by clasper surgery on $G(\b)$.

Finally, let us modify the above discussion for the construction of Theorem 
\ref{thm.2}. Given an $n$-pattern $\b^{(n)}$, let $G(\b^{(n)})$ be a tree
clasper of degree $n$ in $X_0$, which consists of $c^{(n+1)}$ and $2^{n+1}$
leaves $l_i$ (one in each univalent vertex of $c^{(n+1)}$).
There is a 1-1 correspondence between the connected components $T_i$ of 
$\b^{(n)}-c^{(n)}$ and the leaves $l_i$ of $G(\b^{(n)})$. We will choose 
these leaves so that $l_i=\phi(T_i) \in F$, and we will let $L(b^{(n)})$
be obtained from the unlink by clasper surgery on $G(\b^{(n)})$.

We need to show that $L(b^{(n)})$ has trivial $n$th Alexander module. 
Indeed, using the figures above that describe clasper surgery, it follows
that clasper surgery on $G(\b^{(n)})$ is equivalent to surgery on a clasper
$G'(\b^{(n)})$ of degree $1$ whose leaves lie in $F^{(n)}$. This implies
that the $n$th Alexander module of $L(b^{(n)})$ is trivial.

We end this section with a comment on pictures.
To get {\em pictures} of the constructed links, one may use various 
descriptions of surgery on a clasper that were discussed at length by
Goussarov and Habiro at \cite{Gu2,H}. From our point of view though, these
pictures are complicated and unnecessary, since not only claspers describe
surgery adequately, but also the invariants which we will use behave well
with respect to clasper surgery. This is the content of the next section.

\section{Computing the tree part of the Aarhus integral}
\lbl{sec.trees}

\subsection{The Aarhus integral in brief}
\lbl{sub.aarhus}

As was stated in the discussion of Theorem \ref{thm.1}, we will not
compute the Milnor invariants of the links $L(\b)$ constructed via
clasper surgery, but rather we will compute the tree-part of their Aarhus
integral. The Aarhus integral is a graph version of stationary phase
approximation that was introduced at \cite{A}. Despite its intimidating name, 
it is a rather harmless combinatorial object which we now describe.

Consider a framed link $C \subset S^3-\O$ and let $(M,L)=(S^3,\O)_G$ denote 
the result of
surgery on $C$. That is, $M$ is the 3-manifold obtained from $S^3$ by surgery
on $C$ and $L$ is the image of $\O$ after surgery. Assuming that $M$ is a 
rational homology sphere (i.e., that the linking matrix of $C$ has nonzero
determinant) the Aarhus
integral $Z(M,L)$ can be computed by the Kontsevich integral of the link
$\O \cup C$ by integration as follows:
$$
Z(M,L)=\int dX \, Z(S^3,\O \cup C)
$$
(where $X$ is a set of variables in 1-1 correspondence with the components of
$C$).
Let us briefly recall from \cite{A} how this integration works.
Consider an element
$$
s = \exp \left( \frac{1}{2} \sum_{x,y \in X}
 \st{x}{y}{Q_{xy}} \right) R ,
$$
with $R$ a series of graphs that do not contain a strut whose legs are 
colored by $X$. Notice that $Q$ and $R$, the $X$-{\em strutless part} of
$s$, are uniquely determined by $s$. 
Then, the integration $\int \, dX (s)$ glues all the $X$-colored legs
of $R$ pairwise, using the negative inverse of the matrix $Q$. That is,
when two legs $x,y$ of $R$ are glued, the resulting graph is multiplied by
$-Q^{xy}$, the negative inverse of the matrix $Q_{xy}$.

It follows immediately that the {\em tree-part} $\Ztr(M,L)$ of $Z(M,L)$
depends only on the tree-part $\Ztr(S^3,\O \cup C)$ of $Z(S^3,\O \cup C)$.

\subsection{Claspers and the Aarhus integral}
\lbl{sub.link}

Let us adapt the above discussion when the link $C$ is one that describes
clasper surgery. Consider a null clasper $G \sub S^3-\O$ of degree $1$
constructed from a pattern $\b$ and let $(M,L)=(S^3,\O)_{G}$. 
Let $\Zmin(M,L)$ denotes the {\em lowest degree nonvanishing tree part} 
of $\Ztr(M,L)$. Assuming that the pattern is nonvanishing, and after we choose
string-link representatives of $L \cup G$, 
we will show that

\begin{proposition}
\lbl{prop.min}
We have
$$
\Zmin(M,L)=\b \in \At
$$
\end{proposition}

It is clear that this concludes Theorem \ref{thm.1}.

\begin{proof}(of Proposition \ref{prop.min})
Surgery on $G$ is equivalent to surgery on a $6$ component link $C=
C^e \cup C^l$; see Section \ref{sub.surgery}. $C^e$ is a borromean link and 
$C^l$ consists of the leaves of $G$. 
In the obvious basis, the linking matrix of $C$ and its negative 
inverse are given as follows:
$$
\left(\begin{array}{cc}
0 & I \\ I & \lk(C^l_{i},C^l_{j}) \\
\end{array}\right)
\hspace{0.5cm} \text{ and } \hspace{0.5cm}
\left(\begin{array}{cc}
\lk(C^l_{i},C^l_{j}) & -I \\ -I & 0 \\
\end{array}\right) .
$$
In particular, a univalent vertex labeled by a leaf has to be glued to a 
univalent vertex labeled by the corresponding edge.
Let $A_i=\{C^e_{i},C^l_{i}\}$ denote
the arms of $G$ for $i=1,2,3$. It is a key fact that surgery on any proper
subcollection of the set $\{A_1,A_2,A_3\}$ of arms does not change the pair
$(S^3,\O)$. In other words, alternating with respect to the $8$ subsets
of the set of arms we have that
$$
Z([(S^3,\O), G])=Z([(S^3,\O), \{A_1,A_2,A_3\} ])
$$
The nontrivial contributions to the left hand side come from
the $(\O \cup C)$-strutless part of $Z(S^3, \O \cup C)$ 
that consists of graphs with legs on $A_1$ and on $A_2$ and on $A_3$.

What kind of diagrams in $\Ztr(S^3,\O \cup C)$ contribute to the above sum? 
Consider a disjoint union $D$ of trees whose legs are labeled by $\O \cup C$. 
$D$ must have a leg (i.e., univalent vertex) labeled by $C^l_i$ or by $C^e_i$ 
for each $i=1,2,3$. If $D$ has a leg labeled by $C^l_i$, then due to the shape
of the gluing matrix, $D$ must have a $C^e_i$-labeled leg. Thus, in all
cases, $D$ must have legs labeled by all three edges $C^e_i$ of $G$.

Consider a tree $T$ labeled by $\O \cup C$. If $T$ has a $C^e_i$-labeled leg,
then it must either have legs labeled by all three edges of $G$, or else it
must have a leg labeled by $C^l_i$. Indeed, $C^e_i$ is an unknot in a ball
disjoint from $\O \cup C-\{C^l_i\}$, thus the rest of the trees have vanishing
coefficient in $\Ztr(S^3,\O \cup C)$.

Consider further a {\em vortex} $Y$ (that is, a unitrivalent graph of the shape
$Y$ with three univalent vertices and one trivalent one) whose legs are 
labeled by three leaves of $G$. Then, the coefficient of $Y$ in 
$Z(S^3,\O \cup C)$ is $1$.

Consider further a tree $T$ with one univalent vertex labeled by a leaf $C^l_i$
of $G$ and all other vertices labeled by $\O$. Recall the corresponding
rooted tree $T_i$ which is a component of $\b-v$. Then the coefficient of $T$
in $\Ztr(S^3,\O \cup C)$ is zero if $\deg(T) < \deg(T_i)$ and equals to
$1$ if $T=T_i$. This, together with the above discussion and the gluing rules
concludes the proof of Proposition \ref{prop.min}.
The argument is best illustrated by the following figure:
$$
\psdraw{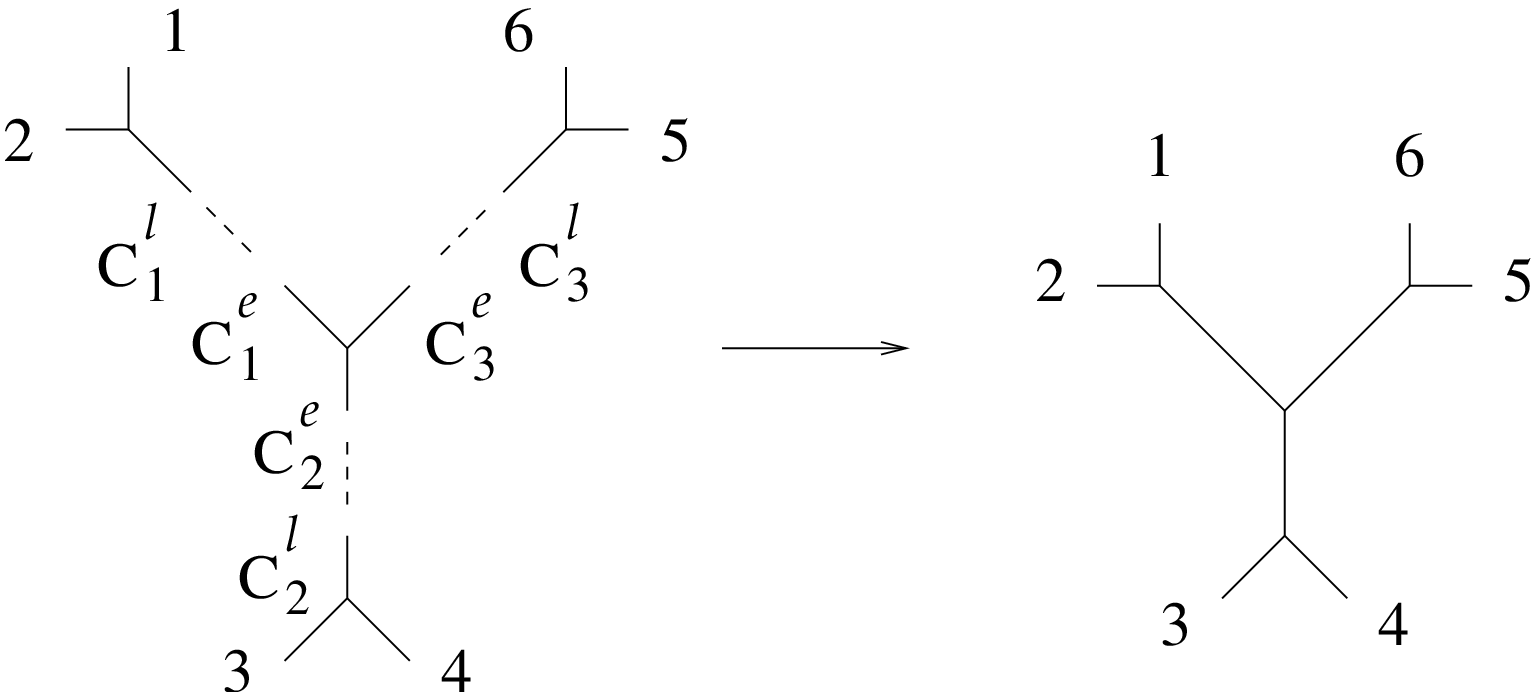}{3in}
$$
%The above discussion generalizes immediately to the case of arbitrary 
%collections of null claspers.
\end{proof}

The above proposition and its proof generalize easily to the case of claspers
$G$ corresponding to nonvanishing $n$-patterns $\b^{(n)}$. In that case,
if $(M,L)$ denote the corresponding link, we still have that
$$
\Zmin(M,L)=\b^{(n)} \in \At
$$
which implies Theorem \ref{thm.2}.

\begin{remark}
In the above discussion we have silently chosen dotted Morse link 
representatives (or equivalently, string-link representatives) and we ought 
to have normalized the Aarhus integral. But this does not affect the lowest 
degree nonvanishing tree part.
\end{remark}

The links constructed by clasper surgery in Theorem \ref{thm.1} include
the links that Cochran constructed via Seifert surfaces.

\begin{question}
Does Section \ref{sub.construct} construct every link with trivial Alexander 
module?
\end{question}

%If we could construct some special Seifert surface for a link with trivial
%Alexander module, then perhaps we could answer the above question 
%positively.

\ifx\undefined\bysame
	\newcommand{\bysame}{\leavevmode\hbox
to3em{\hrulefill}\,}
\fi

\end{document}